\documentclass[a4paper,12pt]{article}
\usepackage{amsmath}
\usepackage{amssymb}
\usepackage{array}

\begin{document}
\title{\bf A property of $C^{k,\alpha}$ functions}
\date{}
\author{{\bf Robert DALMASSO}\\\\ Le Galion - B\^atiment B, \\33 Boulevard Stalingrad, 06300 Nice,  France}

\maketitle
\footnote[0]{\sl E-mail address: {\rm robert.dalmasso51@laposte.net}}

\allowdisplaybreaks[4]

\bigskip

\hrule

\bigskip

\noindent{\small{\bf Abstract.} Let $f$ be a nonnegative function of class $C^k$ ($k \geq 2$) such that $f^{(k)}$ is H\"older continuous  with exponent $\alpha$ in $(0,1]$. If $f'(x) = \cdots = f^{(k)}(x) = 0$ when $f(x) = 0$, we show that $f^{\mu}$ is differentiable for $\mu \in (1/(k+\alpha), 1)$ and under an additional condition we  show that $(f^\mu)'$ is  H\"older continuous with exponent $\beta = \mu(1+\alpha) - 1$ (if $\beta \leq 1$) at $x \in [0,T]$ when $f(x) = 0$. $(f^\mu)'$ is Lipschitz continuous at $x$ if $f(x) > 0$.}

\bigskip



\bigskip

\hrule

\bigskip

\noindent$C^{k}[a,b]$ denotes the space of functions differentiable up to order $k$ such that the derivatives of order $k$ are continuous on $[a,b]$ and $C^{k,\alpha}[a,b]$ denotes the space of functions in $C^{k}[a,b]$ such that the derivatives of order $k$ are H\"older continuous  with exponent $\alpha$ in $(0,1]$. Recall that $g: [a,b] \to \mathbb R$ is H\"older continuous  with exponent $\alpha \in (0,1]$ at $x \in [a,b]$ if
$$\sup\{|g(y)-g(x)||y-x|^{-\alpha}\, ; \, y \not= x \, ,  \, y \, \in \, [a,b]\} < \infty \, ,$$
and that $g$ is H\"older continuous  with exponent $\alpha \in (0,1]$ in $[a,b]$ if
$$\sup\{|g(x)-g(y)||x-y|^{-\alpha}\, ; \, x \not= y \, , \, x \, , \, y \, \in \, [a,b]\}  < \infty \, .$$

It is well-known (\cite{gl})that if a nonnegative function $f$ is in $C^2[a,b]$ and if the second derivative of $f$ vanishes at the zeros of $f$,  then $f^{1/2}$ is in $C^1[a,b]$. Now if $f \in C^m[a,b]$ is nonnegative and if all its derivatives vanish at the zeros of $f$, then $f^{1/m}$ is not necessarily in $C^1[a,b]$ (See \cite{di}). Finally let $f \in C^{k,\alpha}[a,b]$, $k \geq 1$ and $f \geq 0$. Then $f^{1/k+\alpha}$ is absolutely continuous (See \cite{cjs} Lemma 1 and also Remark 2 in \cite{d} when $k = 1$).

\medskip

Now let $f \in C^{k,\alpha}[0,T]$, $T > 0$, $k \geq 2$, be such that $ f^{(j)}(x) = 0$ for some $x \in [0,T]$, $j = 0,\cdots,k$. Then we define

\begin{displaymath}
N(x,y) = \displaystyle (y-x)^{k-1}\int_{0}^{1}(1-s)^{k-2}f^{(k)}(sy+(1-s)x) \, ds \, \, ,
\end{displaymath}
and, if $f \geq 0$,
\begin{displaymath}
D(x,y) = \displaystyle ((y-x)^{k}\int_{0}^{1}(1-s)^{k-1}f^{(k)}(sy+(1-s)x) \, ds)^{(k+\alpha-1)/(k+\alpha)}
\end{displaymath}
for $x \, , \, y \in [0,T]$. We have the following theorem.

\bigskip

\noindent{\bf Theorem.} {\sl Let $f \in C^{k,\alpha}[0,T]$, $T > 0$, $k \geq 2$, be such that $f \geq 0$. Assume that $f$ has at least one zero in $[0,T]$. If $f'(x) = \cdots = f^{(k)}(x) = 0$ when $f(x) = 0$, then $f^{\mu}$ is differentiable for $\mu \in (1/(k+\alpha), 1)$. If moreover $N(x,y)/D(x,y)$ is bounded for $(x,y) \in \{t \in [0,T]; \, f(t) = 0\}\times\{t \in [0,T]; \, f(t) > 0\}$, then $(f^\mu)'$ is  H\"older continuous with exponent $\beta = \mu(k+\alpha) - 1$ at $x$ such that $f(x) = 0$ (if $\beta \leq 1$). $(f^\mu)'$ is Lipschitz continuous at $x$ if $f(x) > 0$.}

\bigskip
\noindent{\sl Proof.} $f^{\mu}$ is clearly differentiable at $x \in [0,T]$ when $f(x) > 0$. Suppose that $f(x) = 0$. For $y \in [0,T]$ we can write
\begin{equation}
\begin{array}{lcl} 
 f(y)   & = & \displaystyle  \frac{(y - x)^k}{(k-1)!}\int_0^1{(1-s)^{k-1}}{f^{(k)}(sy+(1-s)x)\, ds} \\ \\
  &  \leq & \displaystyle \frac{|y - x|^k}{(k-1)!}\int_0^1{(1-s)^{k-1}}{|f^{(k)}(sy+(1-s)x)|\, ds} \\ \\
  & \leq & \displaystyle C\frac{|y - x|^{k+\alpha}}{(k-1)!}\int_0^1{(1-s)^{k-1}}{s^\alpha\, ds} \\ \\
  & = & \displaystyle \frac{C}{(1+\alpha)\cdots(\alpha+k)}|y - x|^{k+\alpha}  \, \, ,
\end{array}
\label{eq:eq1}
\end{equation}
for some constant $C$, which implies that $f^\mu$ is differentiable at $x$.

Let $x \in [0,T]$. Suppose first that $f(x) = 0$. Then $f^{(j)}(x) = 0$ for $j = 1,\cdots,k$. Let $y \in [0,T]$  be such that $f(y) > 0$. We can write
\begin{displaymath}
f'(y) = \displaystyle \frac{(y-x)^{k-1}}{(k-2)!}\int_{0}^{1}(1-s)^{k-2}f^{(k)}(sy + (1-s)x) \, ds \, \, ,
\end{displaymath}
and
\begin{displaymath}
f(y) = \displaystyle \frac{(y-x)^{k}}{(k-1)!}\int_{0}^{1}(1-s)^{k-1}f^{(k)}(sy + (1-s)x) \, ds \, \, .
\end{displaymath}
Using \eqref{eq:eq1} we get
\begin{displaymath}
\begin{array}{lcl} 
 |(f^\mu)'(y) - (f^\mu)'(x)|  & = & \displaystyle  \mu|f(y)^{\mu - 1}f'(y)| \\ \\
 & = & \mu|f(y)^{\mu - \frac{1}{k+\alpha}}f(y)^{-\frac{k+\alpha - 1}{k+\alpha}}f'(y)| \\ \\
 & = & C_1|f(y)^{\mu - \frac{1}{k+\alpha}}|N(x,y)|/D(x,y) \\ \\
  &  \leq & \displaystyle C_{2}f(y)^{\mu - \frac{1}{k+\alpha}} \leq C_{3}|y - x|^\beta \, \, ,\\
\end{array}
\label{eq:eq2}
\end{displaymath}
for some constants $C_{j}$ ($j = 1,\cdots,3$) where $C_2$ and $C_3$ may depend on $x$. Since $f^{\mu}$ is $C^1$ near $t$ when $f(t) > 0$, this implies that $f^{\mu} \in C^1[0,T]$. Suppose now that $f(x) > 0$. There exist $c, d \in [0,T]$ such that $c < d$, $x \in [c,d]$ when $x = 0$ or $x = T$ and $x \in (c,d)$ when $x \in (0,T)$ and $f(y) \geq f(x)/2$ for $y \in [c,d]$. Let $y \in [c,d]$. We have
\begin{displaymath}
\begin{array}{lcl} 
 |(f^\mu)'(y) - (f^\mu)'(x)|  & = & \displaystyle  \mu|f(y)^{\mu - 1}f'(y) - f(x)^{\mu - 1}f'(x)| \\ \\
  &  \leq & \displaystyle  \mu(f(y)^{\mu - 1}|f'(y) - f'(x)| \\ \\
  & & + |f'(x)||f(y)^{\mu - 1} - f(x)^{\mu - 1}|)\\ \\
   & \leq & \displaystyle  C_{1}|y - x| \, \, ,\\
\end{array}
\end{displaymath}
for some constant $C_1$ depending on $x$. Since $(f^\mu)'$ is continuous on $[0,T]$  there exists a constant $C_2$ depending on $x$ such that $|(f^\mu)'(y) - (f^\mu)'(x)| \leq C_{2}|y - x|$ for $y \in [0,T]\backslash[c,d]$.

The proof of the theorem is complete.

\bigskip

\noindent{\bf Remark.} The case $k = 1$ is treated in \cite{d}. Notice that, when $k \geq 2$ and $\mu \in [1/2,1)$, $f^{\mu}$ is in $C^1[0,T]$: See \cite{di} or \cite{gl}. Moreover assume that $k \geq 2$ and that $f'(0) = 0$ (resp. $f'(T) = 0$) when $f(0) = 0$ (resp. $f(T) = 0$). Then, if $\mu \in (1/2,1)$, $(f^\mu)'$ is  H\"older continuous with exponent $2\mu - 1$ at $x$ if $f(x) = 0$ and  Lipschitz continuous at $x$ if $f(x) > 0$: See \cite{d}.
\bigskip

\bigskip

\noindent{\bf Corollary.} {\sl Let $f \in C^{k,\alpha}[0,T]$, $T > 0$, $k \geq 2$. Assume that $f^{(j)}(0) = 0$ for $j = 0,\cdots,k$ and that $f^{(k)} > 0$ on $(0,\eta]$ for some $\eta \in (0,T)$ and $f^{(k)} \geq 0$ on $[\eta,T]$. Then $(f^\mu)'$ is  H\"older continuous with exponent $\beta = \mu(k+\alpha) - 1$  at $0$ (if $\beta \leq 1$). $(f^\mu)'$ is Lipschitz continuous at $x \in (0,T]$.}

\bigskip
\noindent{\sl Proof.} In view of the Theorem it is enough to show that
$N(0,y)/D(0,y)$ is bounded on $(0,T]$. Let
\begin{displaymath}
\displaystyle 0 < \varepsilon < \min(1,(\frac{k-1}{2||f^{(k)}||_{\infty}}\int_0^{1}(1-s)^{k-2}f^{(k)}(sy)ds)^\frac{1}{k-1}) \, \, .
\end{displaymath}
We can write
\begin{displaymath}
\displaystyle \int_0^{1}(1-s)^{k-1}f^{(k)}(sy)ds =  \int_0^{1-\varepsilon}(1-s)^{k-1}f^{(k)}(sy)ds + \int_{1-\varepsilon}^{1}(1-s)^{k-1}f^{(k)}(sy)ds\, \, .
\end{displaymath}
Now we have
\begin{displaymath}
\displaystyle \int_0^{1-\varepsilon}(1-s)^{k-1}f^{(k)}(sy)ds \geq \varepsilon\int_0^{1-\varepsilon}(1-s)^{k-2}f^{(k)}(sy)ds \, ,\end{displaymath}
and
\begin{displaymath}
\displaystyle  \int_{1-\varepsilon}^{1}(1-s)^{k-2}f^{(k)}(sy)ds \leq \frac{\varepsilon^{k-1}||f^{(k)}||_{\infty}}{k-1}\, \, .
\end{displaymath}
Then
\begin{displaymath}
\begin{array}{lcl} 
\displaystyle \int_0^{1}(1-s)^{k-1}f^{(k)}(sy)ds  & \geq & \displaystyle  \varepsilon  \int_0^{1}(1-s)^{k-2}f^{(k)}(sy)ds \\ \\
 & &\displaystyle - \varepsilon \int_{1-\varepsilon}^{1}(1-s)^{k-2}f^{(k)}(sy)ds\\ \\
  &  \geq & \displaystyle\varepsilon  \int_0^{1}(1-s)^{k-2}f^{(k)}(sy)ds - \frac{\varepsilon^{k}}{k-1}||f^{(k)}||_{\infty} \\ \\
  & \geq &\displaystyle \frac{\varepsilon}{2}  \int_0^{1}(1-s)^{k-2}f^{(k)}(sy)ds \, \, .\\
\end{array}
\end{displaymath}
Now, when $y > 0$, we get
\begin{displaymath}
\begin{array}{lcl} 
\displaystyle\frac{N(0,y)}{D(0,y)}  & \leq & \displaystyle  y^{-\frac{\alpha}{k+\alpha}}(\frac{2}{\varepsilon})^{\frac{k+\alpha-1}{k+\alpha}}(\int_0^{1}(1-s)^{k-2}f^{(k)}(sy)ds)^{\frac{1}{k+\alpha}} \\ \\
   &  \leq & \displaystyle C_{1}(\varepsilon)y^{-\frac{\alpha}{k+\alpha}}(y^{\alpha}\int_0^{1}(1-s)^{k-2}s^{\alpha})^{\frac{1}{k+\alpha}}  \leq C_{2}(\varepsilon) \, \, .\\
\end{array}
\end{displaymath}
Then the result follows from the Theorem.

\bigskip

\noindent{\bf Example 1.} Let
\begin{displaymath}
\beta_0 = 0 \, ,\quad \beta_j = \displaystyle \frac{1}{j+1}(\beta_{j-1} + \frac{1}{(j+1)!}) \, \, , j = 1, \cdots ,k  \quad \textrm{and} \quad T \in (0,1]\, \, ,
\end{displaymath}
and let
\[ f(x) = 
\begin{cases}
\displaystyle -\frac{x^{k+1}}{(k+1)!}\ln x + \beta_{k}x^{k+1}& \textrm{if} \quad x \in (0,T] \\
0 & \textrm{if} \quad x = 0\, \, .
\end{cases} \]
Then $f \in C^{k,\alpha}[0,T]$ for all $\alpha \in (0,1)$, $f^{(j)}(0) = 0$ for $j = 0,\cdots, k$ and $f^{(k)}(x) = -x\ln x$. Then we can apply the Corollary. Notice that here  $N(0,y)/D(0,y)$ is continuous on $(0,T]$ and tends to 0 as $y \to 0$.

\bigskip

\noindent{\bf Example 2.} For $\alpha \in (0,1]$ let $f(x) = \displaystyle x^{k+\alpha}g(x)$, $x \in [0,T]$ where $g \in C^{k,\alpha}[0,T]$ is such that $g  > 0$ on $(0,T]$. Then $f \in C^{k,\alpha}[0,T]$, $f^{(j)}(0) = 0$ for $j = 0,\cdots, k$ and $N(0,y)/D(0,y)$ is continuous on $(0,T]$. Suppose that $g^{(j)}(0) \not= 0$ for some $j \in \{0,\cdots,k\}$ and $g^{(i)}(0) = 0$ for $i = 0,\cdots,j-1$ if $j \geq 1$. Then $N(0,y)/D(0,y) \to l$ as $y \to 0$ where $l > 0$ if $j = 0$ and $l = 0$ if $j \in \{1,\cdots,k\}$.

\end{document}